\input amstex
\input xy
\input epsf
\xyoption{all}
\documentstyle{amsppt}
\document
\magnification=1200 
\NoBlackBoxes
\pageheight{17cm}
\hsize14truecm

\def\L{\Cal{L}}
\def\P{\Cal{P}}

\def\Q{\bold{Q}}
\hoffset=0.5in
\voffset=0.5in
\nologo

    

\bigskip
\centerline{\bf COMBINATORIAL CUBIC SURFACES}

\smallskip

\centerline{\bf AND  RECONSTRUCTION THEOREMS}

\bigskip
\centerline{\bf Yu.~I.~Manin}

\medskip

\centerline{\it Max--Planck--Institut f\"ur Mathematik, Bonn, Germany,}

\centerline{\it and Northwestern University, Evanston, USA}

\bigskip

{\bf Abstract.}  This note contains a solution to the following problem:
reconstruct the definition field and the equation of a projective cubic
surface, using only combinatorial information about the set of
its rational points. This information is encoded in two relations:
collinearity and coplanarity of certain subsets  of points. 
We solve this problem, assuming mild ``general position'' properties.

This study is motivated by an attempt to address the Mordell--Weil
problem for cubic surfaces using essentially model theoretic methods.
However, the language of model theory is not used explicitly.

\bigskip

{\bf Contents}

\medskip

0. Introduction and overview

\smallskip

1. Quasigroups and cubic curves

\smallskip

2. Reconstruction of the ground field and a cubic surface
from combinatorics of tangent sections

\smallskip

3. Combinatorial and geometric cubic surfaces

\smallskip

4. Cubic curves and combinatorial cubic curves over large fields

\smallskip

APPENDIX. Mordell--Weil and height: numerical evidence

\smallskip

References

\bigskip

\centerline{\bf 0. Introduction and overview}

\medskip

{\bf 0.1. Cubic hypersurfaces.} Let $K$ be a field, finite or infinite.
In the finite case we assume cardinality of $K$ to be sufficiently
large, the exact lower boundary depending on various particular
combinatorial construction.

\smallskip

Let $P=\bold{P}^N_K$ be a projective space over $K$, with a projective
coordinate system $(z_1:z_2:\dots :z_{N+1})$. A {\it cubic hypersurface $V\subset P$}
defined over $K$ is, by definition, the closed subscheme defined by
an equation $c=0$ where $c\in K[z_1:z_2,\dots :z_{N+1}]$
is a non--zero cubic form. There is a bijection between the set of such subschemes and the set
$\bold{P}^9(K)$ of coefficients of $c$ modulo $K^*$.

\smallskip

We will say that $V$ is  {\it generically reduced} if
after extending $K$ to an algebraic closure $\overline{K}$,
$c$ does not acquire a multiple factor.

\smallskip

In this paper, I will be interested in the following problem: 

\medskip

{\bf 0.1.1. Problem.} {\it Assuming $V$ generically reduced, 
reconstruct  $K$ and the subscheme
$V\subset P= \bold{P}^N_K$  starting with the set of its $K$--points $V(K)$ endowed with
some additional combinatorial structures of geometric origin.} 

\smallskip

The basic combinatorial data that I will be using are subsets of 
smooth points of  $V(K)$ lying upon various sections of $V$ by projective
subspaces of $P$ defined over $K$.  Thus, for the main case treated here, that of
cubic surfaces ($N=3$), I will deal combinatorially with the structure,
consisting of 
\smallskip
{\it a) The subset of smooth (reduced, non--singular) points $S:=V_{sm}(K)$.}

\smallskip

{\it b) A triple symmetric relation ``collinearity'': $\L\subset S^3:=S\times S\times S$.}

\smallskip

{\it c) A set $\P$ of subsets of $S$ called  ``plane sections''.}

\medskip

In the first approximation, one can imagine $\L$ (resp. $\P$) as simply subsets
of collinear  triples (resp.  $K$--points of  $K$--plane sections)  of $V$.
However, various limiting and degenerate cases must be treated with care
as well.
\smallskip

For example, as a working definition of $L$ we  will adopt the following
convention: {\it $(p,q,r)\in S^3$ belongs to $ \L$ if either
$p+q+r$ is the full intersection cycle of $V$ with a $K$--line $l\subset \bold{P}^N$
(with correct multiplicities), or else if there exists a $K$--line $l\subset V$
such that $p,q,r\in l$.}

\medskip

{\bf 0.2. Geometric constraints.}  If an instance of the set--theoretic
combinatorial structure
such as $(S,\L , \P)$ above, comes from  a cubic surface $V$ defined over a field $K$,
we will call such a structure {\it geometric one}.

\smallskip

Geometric structures satisfy additional combinatorial constraints.

\smallskip

{\it The reconstruction problem} in this context consists of two parts:

\medskip

{\it (i) Find a list of constraints ensuring that each   $(S,\L , \P)$ satisfying these constraints
is geometric.} 

\smallskip

{\it (ii) Devise a combinatorial procedure that reconstructs $K$ and
$V\subset \bold{P}^3$ realizing   $(S,\L , \P)$ as a geometric one.}

\smallskip
Besides, ideally we want the reconstruction procedure to be functorial:
certain maps of combinatorial structures, in particular, their isomorphisms,
must induce/be induced by morphisms of ground fields and $K$--linear maps of $P$.

\smallskip

In the  subsection 0.4, I will describe a classical archetype
of reconstruction, -- combinatorial characterization of projective planes.
I will also explain the main motivation for trying to extend this technique
to cubic surfaces: {\it the multidimensional weak Mordell-Weil problem.}

\medskip

{\bf 0.3. Reconstruction of $K$ from curves and configurations of curves.}
One cannot hope to reconstruct the ground field $K$, if $V$
is zero--dimensional or one--dimensional. Only starting with cubic surfaces (N=3),
this prospect becomes realistic.

\smallskip

In fact, if $N=1$, we certainly cannot reconstruct $K$ from any
combinatorial information about one $K$--rational cycle of degree 3 on $P^1_K$.

\smallskip

If $N=2$, then for a smooth cubic curve $V$, the set $V(K)$
endowed with the collinearity relation is the same as $V(K)$
considered as a principal homogeneous space
over the ``Mordell--Weil'' abelian group, unambiguously obtained from $(V(K),L)$
as soon as we arbitrarily choose the identity (or zero) point: cf. a recollection of
classical facts in sec. 1 below.  Generally,
this group does not carry enough information to get hold of $K$, if $K$
is finitely generated over $\bold{Q}$.

\smallskip

However, the situation becomes more promising, if we
assume $V$ geometrically irreducible and having just {\it one
singular point} which is defined over $K$. More specifically,
assume that this point is either an ordinary double point
 with two different branches/tangents
defined over $K$ each, or a cusp with triple tangent line,
which is then automatically defined over $K$. 

\smallskip

In the first case, we will say that $V$ is a curve of {\it multiplicative type}, in the second,
of {\it additive type.}

\smallskip

Then we can reconstruct,
respectively, the multiplicative or the additive group of $K$, up to an isomorphism.
In fact, these two groups are canonically identified with $V_{sm}(K)$
as soon as one smooth $K$--point is chosen, in the same way as the Mordell--Weil group
is geometrically constructed from a smooth cubic curve with collinearity relation.
\smallskip

Finally for $N=3$, now allowing $V$ to be smooth, and under mild genericity restrictions, 
we can combine these two procedures and reconstruct both $K$ and a considerable part
of the whole geometric picture.

\smallskip

The idea, which is the main new contribution of this note,
 is this. Choose two points $(p_m,p_a)$ in $V_{sm}(K)$,
not lying on a line in $V$, whose tangent sections $(C_m,C_a)$ are, respectively,
of multiplicative and additive type.  (To find such points, one might need to replace $K$
by its finite extension first).

\smallskip

Now, one can intersect the tangent planes to $p_m$ and $p_a$ by elements of a 
$K$--rational pencil of planes, consisting of all planes containing $p_m$ and $p_a$.
This produces a birational identification of $C_m$ and $C_a$.

\smallskip

The combinatorial information, used in this construction, can be extracted
from the data $\L$ and $\P$. The resulting combinatorial object,
carrying full information about both $K^*$ and $K^+$, can be then processed into $K$,
if a set of additional combinatorial constraints is satisfied.

\smallskip

Using four tangent plane section in place of two, one can then unambiguously
reconstruct the whole subscheme $V$.

\smallskip

For further information, cf. the main text.

\medskip

{\bf 0.4. Combinatorial projective planes and weak Mordell--Weil problem.}
My main motivation for this study was an analog of Mordell--Weil problem
for cubic surfaces: cf. [M3], [KaM], [Vi].

\smallskip

Roughly speaking, the classical Mordell--Weil Theorem for elliptic curves 
can be stated as follows. Consider a smooth plane cubic curve $C$, i.~e.
a plane model of an elliptic curve, over a field $K$ finitely generated
over its prime subfield. Then the whole set $C(K)$
can be generated in the following way: start with a {\it finite}
subset $U\subset C(K)$ and iteratively enlarge it, adding to already
obtained points  each point $p\circ q\in C(K)$ that is collinear with two points
$p,q\in C(K)$ that were already constructed. If $p=q$, then the third collinear point, by definition,  is
obtained by drawing the tangent line to $C$ at $p$.

\smallskip

In the case of a cubic surface $V$, say, not containing
$K$--lines,  there are two versions of this geometric process
(``drawing secants and tangents'').  We may allow to consecutively add only 
points collinear to $p,q\in V(K)$ when $p\ne q$. Alternatively, 
we may also allow to add all $K$--points of the plane section of
$V$ tangent to $V$ at $p=q$.

\smallskip

I will call the respective two versions of finite generation conjecture
{\it strong, resp. weak, Mordell--Weil problem} for cubic surfaces.

\smallskip

Computer experiments suggest  that weak finite generation might
hold at least for some cubic surfaces defined over $\Q$: see [Vi]
for the latest data. The same experiments indicate however, that
the ``descent'' procedure, by which Mordell--Weil is proved for 
cubic curves, will not work in two--dimensional case: a stable percentage
of $\bold{Q}$--points of height $\le H$ remains not expressible
in the form $p\circ q$ with $p,q$ of smaller height.

\smallskip

In view of this, I suggested in [M3], [KaM] to use a totally different approach to finite
generation, based on the analogy with  classical theory of abstract, or combinatorial, projective planes.

\smallskip

The respective finite generation statement can be stated as follows.

\smallskip
{\it For any field $K$ of finite type over its prime subfield, the whole set
$\P^2(K)$ can be obtained by starting with a finite subset $U\subset  \P^2(K)$
and consecutively adding to it lines through pairs of distinct points,
already obtained, and intersection points of pairs of constructed lines.}

\smallskip

The strategy of proof can be presented as a sequence of the following steps.

\medskip

{\it STEP 1.} Define a {\it combinatorial projective plane} $S,\L$ as 
an abstract set $S$ whose elements are called {\it (combinatorial) points}, endowed with a set of subsets 
of points $\L$ called {\it (combinatorial) lines}, such that each two distinct 
points are contained in a single line, and each two distinct lines
intersect  at a single point.

\medskip

{\it STEP 2.} Find combinatorial conditions upon  $(S,\L)$, that are
satisfied for $K$--points of each geometric projective plane $\P^2(K)$, and that 
exactly characterize geometric planes, so that starting with
$(S,\L)$ satisfying these conditions, one can reconstruct
from $(S,\L)$ a field $K$ and an isomorphism of $(S,\L)$
with $(\P^2(K),\ projective\ K-lines)$ unambiguously. 

\smallskip
In fact, this reconstruction must be also
functorial with respect to embeddings of projective planes
$S\subset S^{\prime}$ and the respective combinatorial lines.

\smallskip

These conditions are furnished by the beautiful {\it Pappus Theorem/Axiom}
(at least, if cardinality of $S$ is infinite or finite but large enough).

\medskip

{\it STEP 3.} Given a geometric projective plane $(\P^2(K),\ projective\ K-lines)$,
start with four points in general position $U_0\subset \P^2(K)$ and generate
the minimal subset $S$ of $\P^2(K)$ stable with respect to drawing lines through two points
and taking intersection point
of two lines. 

\smallskip

This subset, with induced collinearity structure, is a combinatorial projective plane.
It satisfies the Pappus Axiom, because it was satisfied for   $\P^2(K)$.
It is not difficult to deduce then that $S$ is isomorphic to $\P^2(K_0)$,
with $K_0\subset K$ the prime subfield, and the embeddings
$K_0\to K$ and $S\to \P^2(K)$ are compatible with geometry.

\bigskip

{\it STEP 4.} Finally, one can iterate this procedure as follows.
If $K$ is finitely generated, there exists a finite sequence of subfields
$K_0\subset K_1\subset ...\subset K_n=K$ such that each 
$K_i$ is generated over $K_{i-1}$ by one element, say $\theta_i$.
If we already know a finite generating set of points $U_{i-1}
\subset \P^2(K_{i-1})$, define $U_i\subset \P^2(K_i)$
as $U_{i-1}\cup \{(\theta_i:1:0)\}$. One easily sees that  $U_i$
generates $\P^2(K_i)$.

\medskip

{\bf 0.5. Results of this paper.} As was explained in 0.3, results
of this paper give partial versions for cubic surfaces  of Steps 1 and 2 in the 
finite generation proof, sketched above.
I can now reconstruct the ground field $K$ and the total subscheme $V\subset \P^3_K$,
under appropriate genericity
assumptions, from the combinatorics of $V(K)$ geometric origin.

\smallskip

However, these results still fall short of a finite generation statement.

\smallskip

The reader must be aware that this approach is essentially model--theoretic,
and it was inspired by yhe successes of  [HrZ] and [Z].
\smallskip

My playground is much more restricted, and 
I do not use explicitly the (meta)language of model theory,
working in the framework of Bourbaki structures.

\smallskip

More precisely, constructions, explained in sec. 2 and 3, are
oriented to the reconstruction of fields of finite type and cubic
surfaces over them. According to [HrZ] and [Z], if one works over an 
algebraically closed ground field, one can reconstruct combinatorially
(that is, in a model theoretic way) much of the classical algebraic geometry.

\smallskip

In sec. 4, I introduce the notion of a {\it large field}, tailor--made for
cubic (hyper)surfaces, and show that large fields can be reconstructed
even from (sets of rational points of) smooth plane cubic curves,
endowed with collinearity relation and an additional structure consisting of
{\it pencils of collinear points} on such a curve. Any field $K$ having no non--trivial 
extensions of degree 2 and 3 is large, hence large fields lie
between finitely generated and algebraically close ones.

\bigskip
\centerline{\bf 1. Quasigroups and cubic curves}

\medskip

{\bf 1.1. Definition.} {\it Let $S$ be a set and $\L \subset S\times S\times S$ 
be a  subset of triples with the following properties:

\smallskip

(i) $\L$ is invariant with respect to permutations of factors $S$.

\smallskip

(ii) Each pair $p,q\in S$ uniquely determines $r\in S$ such that $(p,q,r)\in \L$.

\smallskip

Then $(S,\L )$ is called a symmetric quasigroup.}
 
\medskip

This structure in fact defines a binary composition law
$$
\circ :\, S\times S\to S:\ p\circ q =r \Longleftrightarrow (p,q,r)\in \L .
\eqno(1.1)
$$
Properties of $\L$  stated in the Definition 1.1 can be equivalently 
rewritten in terms of $\circ$: for all $p,q\in S$
$$
p\circ q = q\circ p,\   p\circ (p\circ q)=q.
\eqno (1.2)
$$
The structure  $(S,\circ )$, satisfying (1.2), will also be called a symmetric quasigroup.
The importance of $\L$ for us is that, together with its versions, it naturally comes from geometry.

\smallskip

In terms of $(S,\circ )$, we can define the following groups. For each
$p\in S$, the map $t_p:\,q\mapsto p\circ q$ is an involutive permutation of $S$:
$t_p^2=\roman{id}_S$.

\smallskip

Denote by $\Gamma =\Gamma (S,\L )$ the group generated by all $t_p, p\in S$.
Let $\Gamma^0\subset \Gamma$ be its subgroup, consisting of products
of an even number of involutions $t_p$.

\medskip

{\bf 1.2. Theorem--Definition.} {\it A symmetric quasigroup $(S,\circ )$
is called abelian, if it satisfies any (and thus all) of the following 
equivalent conditions:

\smallskip

(i) There exists a structure of abelian group on $S$, $(p,q)\mapsto pq$,
and an element $u\in S$ such that for all $p,q\in S$ we have $p\circ q =up^{-1}q^{-1}.$

\smallskip

(ii) The group $\Gamma^0$ is abelian.

\smallskip

(iii) For all $p,q,r\in S$, $(t_pt_qt_r)^2=1$.

\smallskip

(iv) For any element $u\in S$, the composition law $pq := u\circ (p\circ q)$
turns $S$ into an abelian group.

\smallskip

(v) The same as (iv) for some fixed element $u\in S$.

\smallskip

Under these conditions, $S$ is a principal homogeneous space over $\Gamma^0$.}

\medskip

For a proof, cf. [M1], Ch.~I, sec. 1,2, especially Theorem 2.1.

\medskip

{\bf 1.3. Example: plane cubic curves.} Let $K$ be a field, $C\subset \bold{P}^2_K$ an 
absolutely irreducible cubic curve defined over $K$. Denote by $S=C_{sm}(K)\subset C(K)$ the set
of non--singular $K$--points of $C$. Define {\it the collinearity
relation} $\L$ by the following condition: 
$$
(p,q,r)\in \L\ {iff}\  p+q+r\  {is\ the\  intersection\ cycle\ of} \ C\
{with\ a}\ K-{line}.
\eqno(1.3)
$$

\medskip

Then $(S,\L)$ is an abelian symmetric quasigroup. This is a classical result.

\smallskip

More precisely, we have the following alternatives. $C$ might be non--singular 
over an algebraic closure of $K$. Then $C$  is the plane model of an abstract
elliptic curve defined over $K$, the group $\Gamma^0$ can be identified with $K$--points
of its Picard group. We call the latter also the Mordell--Weil group of $C$ over $K$.

\smallskip

Singular curves will be more interesting for us, because they carry
more information about the ground field $K$. Each geometrically irreducible
 singular cubic curve has exactly one singular geometric point, say $p$,
and it is rational over  $K$. More precisely, we will distinguish three cases.

\medskip

(I)  $C$ is  of {\it multiplicative type}. This means that $p$ is a double point two tangents
to which at $p$ are rational over $K$.

\medskip

(II)  $C$ is  of {\it additive type}. This means that $p$ is a cusp: a point with triple tangent.

\medskip

(III)  $C$ is   of {\it twisted type}.  This means that  $p$ is 
a double point $p$ two tangents
to which at $p$ are rational and conjugate over a quadratic extension of $K$.

\medskip

The structure of quasigroups related to singular cubic curves is clarified by the following elementary and well known statement.

\medskip

{\bf 1.3.1. Lemma.} {\it  (i) If $C$ is of multiplicative type,  $\Gamma^0$ is isomorphic to $K^*$.

\smallskip
 (ii) If $C$ is of additive type,  $\Gamma^0$ is isomorphic to $K^+$.
 \smallskip
 (iii) If $C$ is of twisted type,  $\Gamma^0$ is isomorphic to the group of $K$--points
 of a form of $G_m$ or $G_a$ that splits over the respective quadratic extension of $K$. 
 The first case occurs when  $\roman{char} K\ne 2$, the second one when $\roman{char} K=2$.}
\bigskip

{\bf Proof.} {\it (Sketch.)} In all cases, the group law $pq:=u\circ (p\circ q)$, 
for an arbitrary fixed $u\in S$ determines the structure of an algebraic group over $K$
upon the curve $C_0$ which can be defined as the normalization of $C$
with preimage(s) of $p$ deleted. An one--dimensional geometrically 
connected algebraic group becomes
isomorphic to $G_m$ or $G_a$ over any field of definition of its points ``at infinity''.
\smallskip

In the next section, we will recall more precise information about the respective isomorphisms
in the non--twisted cases.

\bigskip
\centerline{\bf 2. Reconstruction of the ground field and a cubic surface}

\smallskip
\centerline{\bf  from combinatorics of tangent sections}

\medskip

{\bf 2.1. The key construction.} Let $K$ be a field of cardinality
$\ge 4.$ Then the set $H:=\bold{P}^1(K)$ consists of $\ge 5$ points.

\smallskip

Consider a family of five pairwise distinct points in $H$ for which we choose the following suggestive notation:
$$
0_a,\infty_a, 0_m, 1_m, \infty_m\in \bold{P}^1 (K).
\eqno(2.1)
$$

In view of its origin, the set $H\setminus \{\infty_a\}$ has a special structure of abelian
group $A$ (written additively, with zero $0_a$). In fact, the choice of any affine coordinate $x_a$
on $\bold{P}^1_K$ with zero at $0_a$ and pole at $\infty_a$ defines this structure:
it sends $p\in H\setminus \{\infty_a\}$ to the value of $x_a$ at $p$, and  addition is addition in $K^+$.
The structure  does not depend on $x_a$, but $x_a$ determines the  isomorphism 
of $G_a$ with $K^+$,
and this isomorphism {\it does} depend on $x_a$: the set of all $x_a$'s
is the principal homogeneous space over $K^*$.

\smallskip

Similarly, the set    $H\setminus \{0_m, \infty_m\}$ has a special structure of abelian
group $M$, with identity $1_m$. A choice of affine coordinate $x_m$
on $\bold{P}^1$, with divisor  supported by$(0_m, \infty_m)$ and taking value $1\in K$ at 
$1_m$,
defines this structure.
Again, it does not depend on $x_m$, but $x_m$ determines its  isomorphism with $K^*$,
and this isomorphism does depend on $x_m$. There are, however, only two choices:
$x_m$ and $x_m^{-1}$. They differ by renaming   $0_m \leftrightarrow \infty_m$.

\smallskip

Having said this, consider now {\it an abstract set} $H$ with a subfamily of five elements
denoted as in (2.1). Moreover,  assume in addition that we are given composition laws
$+$ on $H\setminus \{\infty_a\}$  and $\cdot$ on   $H\setminus \{0_m, \infty_m\}$
turning these sets into two abelian groups, $A$ (written additively, with zero $0_a$) and $M$
(written multiplicatively, with identity $1_m$). Define the inversion map $i:\,M\to M$
using this multiplication law: $i(p)=p^{-1}$.
\smallskip

We will encode this extended version  of (2.1), with  additional data recorded 
in the notation $M,A$, as a bijection
$$
\mu : \ M\cup \{0_m,\infty_m\} \to  A\cup \{\infty_a\}
\eqno(2.2)
$$
It is convenient to extend the multiplication and inversion, resp. addition and sign reversal,
to commutative partial composition laws on two sets (2.2) by the usual rules: for $p\in M$, $q\in A$,
we set
$$
p\cdot 0_m:=0_m,\   p\cdot \infty_m:=\infty_m,\  i(0_m):=\infty_m, \  i(\infty_m):=0_m,
\eqno(2.3)
$$
$$
q\pm \infty_a := \infty_a .
\eqno(2.4)
$$

The following two lemmas  are our main tool in this section.

\medskip

{\bf  2.2. Lemma.} {\it   If (2.2) comes from a projective line as above, then 
the map
$$
\nu :\ M\cup \{0_m,\infty_m\} \to  A\cup \{\infty_a\},
$$
$$
\nu (p):= \mu\{\mu^{-1}[\mu(p)-\mu(0_m)]\cdot i\circ\mu^{-1}[\mu(p)-\mu(\infty_m)]\}
\eqno(2.5)
$$
is a well defined bijection.
\smallskip

Moreover, 
$$
\nu (0_m)= 0_a :=0,\ \nu (\infty_m)= \infty_a :=\infty .
\eqno(2.6)
$$

Finally, identifying $M\cup \{0_m,\infty_m\}$ and   $A\cup \{\infty_a\}$
with the help of $\nu$ and combining addition and multiplication, now (partially) defined
on $H$,
we get upon $H\setminus \{\infty\}$ a structure of the commutative field,
with zero $0$ and identity $1:=\nu (1_m)$. This field is isomorphic to the initial field $K$ .}

\medskip

{\bf Proof.} In the situation (2.1), if $A$ is identified with $K^+$ using an affine coordinate $x_a$, and $M$
is identified with $K^*$ using another affine coordinate $x_m$ as above,
these coordinates are connected by the evident fractional linear transformation,
bijective on $\bold{P}^1(K)$:
$$
x_a = c\cdot (x_m-x_m(0_m))\cdot (x_m-x_m(\infty_m))^{-1},\ c\in K^*.
$$
The definition (2.5) is just a fancy way to render this relation, taking into account that now we
have to add and to multiply in two different locations, passing  back and forth
via $\mu$ and $\mu^{-1}$. Instead of multiplying by $c$, we normalize multiplication so that 
$\nu (1_{\infty})$ becomes identity.

\smallskip

This observation makes all the statements evident.

\medskip

The same arguments read in reverse direction establish the following result:

\medskip
{\bf 2.3. Lemma on Reconstruction.} {\it Conversely, let $M$ and $A$ be two abstract abelian groups,
extended by ``improper elements'' to the sets with partial composition laws 
$M\cup \{0_m,\infty_m\}$ and  $A\cup \{\infty_a\}$, as in (2.3), (2.4).
Assume that we are given a bijection $\mu$ as in (2.2), mapping
$1,0_m$, and $\infty_m$ to $A$. Assume moreover that: 

\smallskip
(i) The respective mapping $\nu$ defined by (2.5) is a well defined bijection.
\smallskip

(ii) The set $A$ endowed with its own addition, and multiplication
transported by $\nu$ from $M$, is a commutative field $K$.
 
\smallskip

Then  we get a natural identification $H=\bold{P}^1(K)$. This construction
is inverse to the one described in sec. 2.1.}

\medskip

{\bf 2.4. Combinatorial projective lines and functoriality.} Let us call
an instance of the data (2.2)--(2.4), satisfying the constraints of Lemma 2.2,
{\it a combinatorial projective line} (this name will be better justified
in the remainder of this section). Let us call triples  $(K, \bold{P}^1(K), j)$
where $j$ is a subfamily of five points in $\bold{P}^1(K)$ as in (2.1),
{\it geometric projective lines.} 

\smallskip

The constructions we sketched above are obviously functorial with respect to
various natural maps such as:

\smallskip

a) {\it On the geometric side:} Morphisms of fields, naturally extended to
projective lines with marked points.  Fractional linear transformations of $\bold{P}^1(K)$,
naturally acting upon $j$ and identical on $K$.

\smallskip

b) {\it On the combinatorial side:} Embeddings of groups $M\to M^{\prime}, A\to A^{\prime}$,
compatible with $(\mu , \mu^{\prime})$ and  on improper points. Automorphisms of  $(M,A)$, 
supplied with compatibly changed $\mu$ and improper points.

\smallskip

These statements can be made precise and stated as equivalence of 
categories. We omit details here. 

\medskip

Now we turn to the description of a bare--bones geometric situation,
that can be obtained (in many ways) from a cubic surface,
directly producing combinatorial projective lines.
\medskip

{\bf 2.5. $(C_m,C_a)$--configurations}. Consider a family of subschemes in $\bold{P}^3_K$,
that we will call {\it a configuration:} 
$$
Conf:=(p_m, p_a; C_m, C_a; P_m, P_a)
\eqno(2.7)
$$ 

It consists of the following data:

\smallskip

{\it (i) Two distinct $K$--points $p_m, p_a\in \bold{P}^3(K).$

\smallskip

(ii) Two distinct $K$--planes $P_m,P_a\subset \bold{P}^3$
such that $p_m\in P_m, p_m\notin P_a$ and 
$p_a\in P_a, p_a\notin P_m.$

\smallskip

(iii) Two geometrically irreducible cubic $K$--curves $C_m\subset P_m$,
$C_a\subset P_a$.}

\medskip
We impose on these data  the following constraints:

\medskip

(A)  $p_m\in C_m(K)$ is a double point, and $C_m$ if of multiplicative type,
in the sense of 1.3.
\smallskip

(B) $p_a\in C_a(K)$ is a cusp, and $C_a$ is of additive type. 

\smallskip

(C) Let $l:=P_m\cap P_a$. Denote by $0_m, \infty_m \in l$
the intersection points with $l$ of two tangents to $C_m$ at $x_m$
(in the chosen order).
Denote by $0_a\in l$  the intersection point with $l$ of the tangent to $C_a$
at $x_a$. These three points are pairwise distinct.
 
\medskip

Let $M:=C_{m,sm}(K)$, $A:=C_{a,sm}(K)$ be the respective sets of
smooth points, with their group structure, induced by collinearity
relation and a choice of $1_m$, resp. $0_a$, as in sec. 1.

\smallskip
Define the bijection $\alpha :\,\widetilde{C}_m(K) \to l(K)$,
where $\widetilde{C}_m$ is the normalization of $C_m$, by
mapping each smooth point  $q\in C(K)$ to the intersection point
with $l$ of the line, passing through $p_m$ and $q$.
The two tangent lines at $p_m$ define the images of two points 
of $\tilde{C}_m$ lying over $p_m$.

\smallskip

Similarly, define the bijection $\beta :\,C_a(K) \to l(K)$,
by mapping each smooth point  $q\in C(K)$ to the intersection point
with $l$ of the line, passing through $p_a$ and $q$. The point on $l$
where the triple tangent at cusp intersects it, is denoted $\infty_a$. 

\smallskip
Finally, put
$$
\mu := \beta^{-1}\circ \alpha :\  \ M\cup \{0_m,\infty_m\} \to  A\cup \{\infty_a\}
\eqno(2.8)
$$

\smallskip

Thus $l(K)$ acquires both structures: of a combinatorial  line and of a geometric
line.

\medskip

{\bf 2.6. $(C_m,C_a)$--configurations from cubic surfaces.} Let $V$ be a smooth
cubic surface defined over $K$. At each non--singular point  $p\in V(K)$, there exists
a well defined tangent plane to $V$ defined over $K$.
The intersection of this plane with $V$, for $p$ outside of a proper Zariski closed subset,
is a geometrically irreducible curve $C$, having $p$ as its single singular point.

\smallskip

Again, generically it is of twisted multiplicative type, if $\roman{char} K\ne 2$,
and of twisted additive type, when $p$ lies on a curve in $V$.

\smallskip

Therefore, under these {\it genericity conditions}, replacing $K$ by its finite extension if need be, and renaming this new field $K$,
we can find two tangent plane sections of $V$ that form a $(C_m,C_a)$--configuration
in the ambient projective space. 

\medskip

{\bf 2.6.1. Example.} Consider the diagonal cubic surface
$\sum_{i=1}^4 a_iz_i^3=0$ over a field $K$ of characteristic $\ne 3$.
Then the discriminant of the quadratic equation defining directions
of two tangents of the tangent section at $(z_1:z_2:z_3:z_4)$,
up to a factor in $K^{*2}$, is
$$
D:= \prod_{i=1}^4 a_iz_i.
$$
Hence the set of points of (twisted) additive type consists
of four elliptic curves
$$
E_i:\  z_i=\sum_{j\ne i}a_jz_j^3=0, \quad i=1,\dots ,4.
$$ 
The remaining  points (outside 27 lines) are of (twisted) multiplicative type. 
Those for which $D\in K^{*2}$ are of purely multiplicative type.

\medskip

{\bf 2.7. Reconstruction of the configuration itself.} Returning to the map (2.8),
we see that $K$ can be reconstructed from the $(C_m,C_a)$ configuration,
using {\it only} the collinearity relation on the set
$$
\widetilde{C}_m(K)\cup C_a(K)\cup l(K).
\eqno(2.9)
$$
Moreover, we get the canonical structure of a projective line over $K$ on $l$,
together with the family of five $K$--points on it.

\medskip

To reconstruct
the whole configuration, as a $K$--scheme up to an isomorphism, 
from the same data, it remains
to give in addition two $0$--cycles on $l$: its intersection
with  $C_m$ and $C_a$ respectively. Again, passing to a finite extension of $K$,
if need be, we may and will assume that  all intersection points in 
$C_m\cap l$, $C_a\cap l$ are defined over $K$. This again means that
these cycles belong to the respective collinearity relation on
$C_m(K)\cup C_a(K)\cup l(K)$.
 
 \smallskip
 
To show that knowing these cycles, we can reconstruct $C_m$ and $C_a$
in their respective projective planes, let us look at the equations
of these curves. 

\smallskip

In $P_m$, choose projective coordinates $(z_1:z_2:z_3)$ over $K$ in such a way
that $l$ is given by the equation $z_3=0$,  $p_m$ is $(0:0:1)$,  equations
of two tangents at $p_m$ are $z_1=0$, $z_2=0$, and the points $0_m,\infty_m$
are respectively $(0:1:0)$ and $(1:0:0)$
Then the equation of $C_m$ must be of the form
$$
z_1z_2z_3 +c(z_1,z_2)=0
$$
where $c$ is a cubic form. To give the intersection $C_m\cap l$ is the same as to give
the linear factors of $c$. Since $z_i$ are defined up to multiplication
by constants from $K^*$, this defines $(C_m,P_m)$ up to isomorphism.

\smallskip

Similar arguments work for $C_a$; its equation in coordinates 
$(z_1^{\prime}:z_2^{\prime}:z_3^{\prime})$ on $P_a$  such that $l$ is defined
by $z_3^{\prime}=0$, will now
be
$$
z_1^{\prime 2}z_3 +c^{\prime}(z_1^{\prime}  ,z^{\prime}_2)=0.
$$
We may normalize $z_2^{\prime}$ by the condition that $0_a=(1:0:0)$,
and then reconstruct linear factors of $c^{\prime}$ from the respective
intersection cycle $C_a\cap l$.

\medskip

{\bf 2.8. Reconstruction of $V$ from a tangent tetrahedral configuration.}
Let now $V$ be a cubic surface over $K$. Assume 
that $V(K)$ contains four points $p_i, i=1,\dots ,4$, such that tangent plains $P_i$ at them
are pairwise distinct. Moreover, assume that  tangent sections
$C_i$ are   either of multiplicative, or of additive type,
and each of these two types is represented by some $C_i$.
One can certainly find such $p_i$ defined over a finite extension of $K$.

\smallskip

We will call such a family of subschemes  $(p_i,C_i, P_i)$ {\it a tetrahedral configuration},
even when we do not assumed a priori that it comes from a $V$. If it comes from a $V$,
we will say that it is {\it a tangent tetrahedral
configuration.}

\smallskip

Without restricting generality, we may choose in the ambient $\bold{P}^3_K$ a coordinate system $(z_1:\dots :z_4)$
 in such a way that $z_i=0$ is an equation of
$P_i$.

\smallskip

If the configuration is tangent to $V$, let $F(z_1,\dots ,z_4)=0$ be the equation of  $V$.
Here $F$ is a cubic
form with coefficients in $K$ determined by $V$ up to a scalar factor.
For each $i\in \{1,\dots ,4\}$, write $F$ in the form
$$
F= \sum_{a=0}^3 z_i^af^{(i)}_{3-a} (z_j\,|\,j\ne i) ,
\eqno(2.10)
$$
where $f_b^{(i)}$ is a form of degree $b$ in remaining variables.

\smallskip

\smallskip

Clearly, $f_3^{(i)}=0$ is an equation of $C_i$ in the plane $P_i$. Hence $K$ and 
this equation 
can be reconstructed, up to a common factor, from a part of the tetrahedral configuration
consisting of $P_i$, another plane $P_j$ with tangent section of different type,
and the induced relation of collinearity on them.

\smallskip

Consider the graph $G=G(V;p_1,\dots ,p_4)$ with four vertices labeled $(1,\dots ,4)$, in which
$i$ and $j\ne i$ are connected by an edge, if there is a cubic monomial
in $(z_k\,|\,k\ne i,j)$, that enters with nonzero coefficients in both $f_3^{(i)}$ and $f_3^{(j)}.$
We want this graph to be connected. This will hold, for example, if
in $F$ all four coefficients at $z_i^3$ do not vanish. 
It is clear from this remark that connectedness of $G$ is an open condition
holding on a Zariski dense subset of all tangent configurations.

\medskip

{\bf 2.8.1. Proposition.} {\it If the tetrahedral configuration is tangent to $V$,
with connected graph $G$, then this $V$ is unique.}

\smallskip

{\bf Proof.} Let $g^{(i)}$ be a cubic form in $z_k, k\ne i,$ such that
$z_i=0, g^{(i)}=0$ are equations of $C_i$. We may change $g^{(i)}$
multiplying them by non--vanishing constants $c_i\in K$. If our configuration is tangent to $V$,
given by (2.10),
we may find  $c_i$ in such a way that $c_ig^{(i)}=f_3^{(i)}$.
The obtained   family of forms  $\{c_ig^{(i)}\}$   is {\it compatible} in the following
sense: if a cubic monomial in only two variables has non--zero coefficients
in two $g^{(i)}$'s, then these coefficients coincide. In fact, they are equal
to the coefficient of the respective monomial in $F$. 

\smallskip

Conversely, if such a compatible system exists, and moreover, the graph
$G$ is connected, then $(c_i)$ is unique up to a common factor.
From such $c_ig^{(i)}$ one can reconstruct a cubic form
of four variables, which will be necessarily proportional to $F$:  coefficient at any cubic monomial $m$ in  
$(z_1,\dots ,z_4)$ in it will be equal to the coefficient
of this monomial in any of  $c_ig^{(i)}$, for which $z_i$ does not divide $m$. 

\medskip

{\bf 2.9. Summary.}  This section was dedicated to several key constructions
that show how and under what  conditions a cubic surface $V$ considered as a scheme,
together with a ground field $K$, can be reconstructed
from its set of $K$--points, endowed with some combinatorial data.

\smallskip

The main part of the data was the collinearity relation on $V_{sm}(K)$,
and this relation, when it came from geometry, satisfied some strong
conditions stated in Lemma 2.2.

\smallskip

However, this Lemma and the data used in 2.8 made appeal
also to information about points on the lines of intersections of tangent planes:
cf. specifically constructions of maps $\alpha$ and $\beta$
before  formula (2.8).

\smallskip

We want to get rid of this extra datum and work {\it only} with  points
of $V_{sm}(K)$. 

\smallskip

This must be compensated by taking in account, besides the
collinearity relation, an additional {\it coplanarity relation}
on $V(K)$, essentially given by the sets of $K$--points of
(many) non--tangent plane sections.

\smallskip

The next section is dedicated principally to  a description of the relevant abstract
combinatorial framework. The geometric situations are used
mainly to motivate or illustrate  combinatorial
definitions and axioms.

\bigskip

\centerline{\bf 3. Combinatorial and geometric cubic surfaces}

\medskip

{\bf 3.1. Definition.} {\it A combinatorial cubic surface is an abstract set $S$
endowed with two structures: 
\smallskip
(i) A symmetric  ternary relation  ``collinearity'': $\L\subset S^3 $.
We will say that triples $(p,q,r)\in \L$ are collinear.

\smallskip

(ii) A set $\P$ of subsets  $C\subset S$ called plane sections.}

\smallskip
These relations must satisfy the axioms made explicit below in the subsections 3.2 and  3.3.
Until all the axioms are stated and imposed, we may call a structure
$(S,\L,\P)$ a {\it cubic pre--surface.}

\medskip

{\bf 3.2. Collinearity Axioms.} {\it (i) For any $(p,q)\in S^2$, there exists 
an $r\in S$ such that $(p,q,r)\in L$.

 Call the triple $(p,q,r)$ strictly collinear, if  $r$ is unique  with this property,
and $p,q,r$ are pairwise distinct.

\smallskip
(ii) The subset $\L_s\subset \L$ of strictly collinear triples is 
a symmetric ternary relation.

\smallskip

(iii) Assume that $p\ne q$ and that there are two distinct $r_1,r_2\in S$
with $(p,q,r_1)\in \L$ and $(p,q,r_2)\in \L$. Denote by $l=l(p,q)$ the set of all such
$r$'s. Then $l^3\subset \L$, that is any  triple $(r_1,r_2,r_3)$ of points
in $l$ is collinear.

\smallskip

Such sets $l$ are called lines in $S$.}

\medskip

{\bf 3.2.1.  Example: combinatorial cubic surfaces of geometric origin.} Let $K$ be a field, 
and $V$  a cubic surface
in $\bold{P}^3$ over $K$. Denote by $S=V_{sm}(K)$ the set of nonsingular $K$--points of $V$.

\smallskip 
We endow $S$ with the following relations:

\smallskip

(a) $(p,q,r)\in \L$ iff either $p+q+r$ is the complete intersection cycle of $V$ with a line
in $\bold{P}^3$ defined over $K$ ($K$--line), or else if $p,q,r$ lie on a $K$--line $\bold{P}^1_K$,
entirely contained in $V$.

\smallskip

(b) Let $\bold{P}\subset \bold{P}^3$ be a $K$--plane. Assume that it either contains at least two distinct points of $S$, or is tangent to a $K$--point $p$, or else contains the tangent line
to one of the branches of the tangent section of multiplicative type.
Then 
$C:=\bold{P}(K)\cap S$ is an element of $\P$. All elements of $\P$ are obtained in this way.

\medskip

{\bf 3.3. Plane sections.} We now return to the general 
combinatorial situation. Let $(S,\L,\P)$
be a cubic pre--surface.

\smallskip
 
For any $p\in S$, put
$$
C_p=C_p(S):= \{ \,q\,|\, (p,p,q)\in \L\} \cup \{p\}.
\eqno (3.1)
$$
\medskip

{\bf 3.3.1. Tangent Plane Sections Axiom.} {\it For each $p\in S$,
we have $C_p\in \P.$ Such plane sections are called tangent ones.}

\medskip

The next geometric property of plane sections of geometric cubics
can now be rephrased combinatorially as follows.
\medskip

{\bf 3.3.2.  Composition Axiom.}     {\it (i)  Let $C\in \P$
be a non--tangent plane section containing no lines in $S$.
Then the collinearity relation $\L$ induces on such $C$
a structure of Abelian symmetric quasigroup (cf. Theorem--Definition 1.2).

\smallskip

(ii) Let $C_p=C_p(S)$ be a tangent plane section containing no lines.
Then $\L$ induces on 
$C^0_p:=C_p\setminus \{p\}$ a structure of Abelian symmetric quasigroup.}

\medskip

Choosing a zero/identity point in $C$, resp. $C_p\setminus \{p\}$, we get in this way
a structure of abelian group on each of these sets.

\medskip

{\bf 3.3.3. Pencils of Plane Sections Axiom.} {\it Let $\lambda := (p,q,r)\in \L$. Assume that at least two
of the points $p,q,r$ are distinct. Denote by $\Pi_{\lambda}\subset \P$ the set
$$
\Pi_{\lambda}:= \{C\in \P\,|\, p,q,r\in P\}.
\eqno(3.2)
$$
and call such $\Pi_{\lambda}$'s pencils of plane sections.

Then we have:

\smallskip

(i) If $(p,q,r)$ do not lie on a line in $S$, then
$$
S\setminus \{p,q,r\}=\coprod_{C\in \Pi_{\lambda}} (C\setminus \{p,q,r\})
\eqno(3.3)
$$
(disjoint union).

\smallskip

(ii)   If $(p,q,r)$   lie on a line $l$, then
$$
S\setminus l =\coprod_{C\in \Pi_{\lambda}} (C\setminus l )
\eqno(3.4)
$$ 
(disjoint union).}

\medskip

{\bf 3.4. Combinatorial plane sections $C_p$ of multiplicative/additive types.} 
First of all we must postulate $(p,p,p)\in \L$, since in the geometric case $(p,p,p)\notin \L$
can happen only in a twisted case.

\smallskip

There are two different approaches to the tentative distinction between multiplicative
and additive types.
In one, we may try to 
prefigure the future realization of $C_m$ and $C_a$ as essentially the 
multiplicative (resp. additive) groups of a field $K$ to be constructed.

\smallskip
Then, restricting ourselves for simplicity by fields of characteristic zero, we see
that $C_p\setminus \{p\}$ which is of additive type after a choice of $0_a$
must become a vector space over $\bold{Q}$  (be
uniquely divisible), whereas the respective group of multiplicative
type is never uniquely divisible.

\smallskip

However, these restrictions are too weak.

\smallskip

Instead, we will define {\it pairs} of combinatorial tangent plane sections
modeled on $(C_m,C_a)$--configurations of sec. 2. After this is done,
we will be able to ``objectively'', independently of another
member of the pair,  distinguish between $C_m$ and $C_a$
using e.g. the divisibility criterion.

\medskip

{\bf 3.5. Combinatorial $(C_m,C_a)$--configurations.} 
We can now give a combinatorial version of those $(C_m,C_a)$--configurations,
that in the geometric case consist of two tangent plane sections of a cubic surface,
one of additive, another of multiplicative type.

\smallskip

The main point is to see, how to use combinatorial plane sections
in place of ``external'' lines $l=P_m\cap P_a$. This is possible,
because the set of points of this line will now be replaced by bijective to it
set of plane sections, belonging to a pencil, defined in terms of $(C_m,C_a)$,
and geometrically consisting just of all sections containing $p_m$ and $p_a$.

\medskip

Let $(S,\L,\P)$ be a combinatorial pre--surface, satisfying 
Axioms 3.2, 3.3.1, 3.3.2, 3.3.3.

\smallskip

Start with two distinct points of $S$,
not lying on a line in $S$, and respective tangent sections of $S$
$$
(p_m,p_a; C_{p_m},C_{p_a})
\eqno(3.5)
$$

\smallskip
Let $r\in S$ be the unique third point such that
$(p_m,p_a,r)\in \L$, $\lambda := \{p_m,p_a,r\}.$ Put
$C^0_{p_m}:=C_{p_m}\setminus \{p_m\}$,
$C^0_{p_a}:=C_{p_a}\setminus \{p_a\}$.

\smallskip

Denote by $\Pi_{\lambda}$ the respective pencil of plane sections.
Consider the following binary relation:
$$
R\subset C_{p_a}\times C_{p_m}:\ (p,q)\in R \Longleftrightarrow
\exists P\in \Pi_{\lambda},\, p,q\in P.
\eqno(3.6)
$$

\medskip

{\bf 3.5.1. Definition.} {\it $(p_m,p_a; C_{p_m},C_{p_a})$ is called a
$(C_m,C_a)$--configuration, if
the following conditions are satisfied.

\smallskip

(i) $R$ is a graph of some function
$$
\lambda :\ C_{p_a}\to C_{p_m}
\eqno(3.7)
$$
This function must be a bijection outside of two  distinct points
$0_m,\infty_m\in C_{p_a}^0$ which are mapped to $p_m$. Besides,
we must have $\lambda (p_a)\in C^0_{p_m}$.

\smallskip

Assuming (i), put
$$
A:=C^0_{p_a},\ M:= C^0_{p_m}
$$
Introduce on these sets the structures of abelian groups
using the  the Composition Axiom 3.3.2 and some choices of zero and identity
$$
0_a\in C^0_{p_a},\ 1_m\in C^0_{p_m}.
$$
Define the map
$$
\mu : \ M\cup \{0_m,\infty_m\} \to  A\cup \{\infty_a\}
\eqno(3.8)
$$
which is $\lambda^{-1}$ on $M$ and identical on  $0_m,\infty_m$.   
Then 

\smallskip
(ii) Conditions of Lemma 2.2 must be satisfied for this $\mu$.
\smallskip

(iii) $C_{p_m}\cap C_{p_a}$ consists of three pairwise distinct points.}

\medskip

Thus, if  $(p_m,p_a; C_{p_m},C_{p_a})$ is  a
$(C_m,C_a)$--configuration, then
  we can combinatorially reconstruct the ground field and the isomorphic
geometric configuration.

\smallskip

However, passing to tetrahedral configurations, we have to impose 
additional  combinatorial compatibility
conditions, that in the geometric case were automatic.

\smallskip

They are of two types:

\medskip

(a) If two planes $P_i$, $P_j$ of the tetrahedron carry plane sections
of the same type (both additive or both multiplicative), we must
write combinatorially maps, establishing their isomorphism, and postulate this 
fact in the combinatorial setup.
\smallskip
This can be done  similarly to the case of
$(C_m,C_a)$--configurations.  
\medskip
(b) If a schematic tangent plane section $C_i$ can be reconstructed from
two different pairs of tetrahedral plane sections $C_i, C_j$
and $(C_i,C_k)$, the results must be naturally isomorphic.  
\medskip
It is clear, how to do it in principle, and the respective constraints
must be stated explicitly. 
\smallskip

I abstain from elaborating all details here for the following  reason.

\smallskip
If, as a main application of this technique, one tries to imitate the approach to weak Mordell--Weil problem
following the scheme of 0.4, then the necessary combinatorial constraints will
probably hold automatically for finitely generated
combinatorial subsurfaces of an initial geometric
surface. 
\smallskip
The real problem is: how to recognize that a given (say, finitely generated)
combinatorial subsurface of a geometric surface is actually
the whole geometric surface.

\smallskip

I do not know any answer to this problem. 

\smallskip

It is well known, however, that such proper combinatorial 
subsurfaces do exist. 

\smallskip

For example, when $K=\bold{R}$ and $V(\bold{R})$ is not connected, one of the components
can be a combinatorial cubic surface in its own right. More generally,
some unions of classes of the universal equivalence relation ([M1])
are closed with respect to collinearity and coplanarity relations:
this can be extracted from  the results of [SwD1].

\bigskip

\centerline{\bf 4. Cubic curves and combinatorial cubic curves over large fields}

\medskip

{\bf 4.1. Large fields and smooth cubic curves.}  Consider a smooth cubic curve
$C\subset \P^2_K$ defined over $K$. Put  $S:=C(K)$ and endow $S$ with the
collinearity relation $\L\subset S^3$ defined by (1.3). Let $\L^0:=L/\bold{S}_3$,
the set of orbits of $\L$ with respect to the permutations. We may and will represent
the image in $\L^0$ of $(p,q,r)\in \L$ as the 0--cycle $p+q+r$. 

\smallskip
Now assume that $K$ has no non--trivial extensions of degree 2 and 3.
Then  
all intersection points of any $K$--line with $C$ lie in $C(K).$
Therefore,  we have the canonical bijection
$$
\xi :\ \{K-lines\  in\ {\bold{P}}^2(K)\} \to \L^0:\quad {l} \mapsto\ intersection\ cycle\  l\cap C.
\eqno(4.1)
$$
In this approach, $K$--points of  ${\bold{P}}^2(K)\}$ have to be characterized
in terms of pencils of all lines passing through a given point. Therefore, it is more natural to work with the dual projective plane from the start.

\smallskip

Let $\widehat{\bold{P}}^2$ be the projective plane dual to the plane
in which $C$ lies. Combinatorially, $K$--points $\widehat{l}$  (resp. lines
$\widehat{p}$) of  $\widehat{\bold{P}}^2$
are $K$--lines $l$ (resp. points $p$) of  $\P^2_K$, with inverted incidence relation:
$\widehat{l}\in  \widehat{p}$ iff $p\in l$. Thus, (4.1) turns into a bijection
$$
\widehat{\xi }:\ \{K-points\  in\ \widehat{\bold{P}}^2(K)\} \to \L^0:\quad \widehat{l} \mapsto\ intersection\ cycle\  l\cap C.
\eqno(4.2)
$$
Thus, $\widehat{\xi}$ sends   lines in   $\widehat{\bold{P}}^2$ to certain subsets in $\L^0$
that we also may call {\it pencils.} 

\smallskip

This geometric situation motivates the following definition.

\smallskip

Let $(S,\L )$   be an abelian symmetric quasigroup. Put  $\L^0=\L/\bold{S}_3$.

\smallskip
Assume that $\L^0$ is endowed  with a set of its subsets  $\P^0$, called {\it pencils, }
which turns it into a combinatorial projective plane, with pencils as lines.
This means that besides the trivial incidence conditions,
Pappus Axiom is also valid. Hence we can reconstruct from $(\L^0, \P^0)$
a field $K$ such that $\L^0=\bold{P}^2(K)$, $\P^0 =$ the set of $K$--lines
in   $\bold{P}^2(K)$. 

\smallskip

The following Definition, inspired by geometry, encodes the  interaction
between the structures $(S,\L)$ and   $(\L^0, \P^0)$.

\medskip

{\bf 4.2. Definition.} {\it The structure $(S,\L,\P^0 )$ is called a combinatorial
cubic curve over a large field, if the following conditions are satisfied. 

\smallskip

(i) For each fixed $p\in S$, the set of cycles $p+q+r \in \L^0$, $q,r\in S$, is a pencil $\Pi_p$.

\smallskip

(ii) If a pencil $\Pi$  is not of the type $\Pi_p$, then any two distinct elements
in $\Pi$ do not intersect (as unordered triples of $S$--points).

\smallskip

(iii) For each pencil $\Pi_q$ (resp. each pencil not of type $\Pi_q$) and any $p\in S$,
(resp. any $p\ne q$) there exists a unique
cycle in $\L^0$ contained in $\Pi$ and containing $p$.}

\medskip

Obviously, each geometric smooth cubic curve over a large field 
defines the respective combinatorial object.

\medskip

{\bf 4.3. Question.}  Are there such combinatorial curves not coming from
a geometric one? In particular, are fields $K$ coming from 
such combinatorial objects necessarily ``large'' in the sense of algebraic
definition above (closed under taking square and cubic roots)?

\medskip

Similar constructions can be done and question asked for cubic surfaces.
Notice that over a large field, any point on a smooth surface, not lying on a line,
is of either multiplicative, or additive type.

\bigskip

\centerline{\bf APPENDIX. Mordell--Weil and height: numerical evidence}

\medskip

Let $V$ be a geometrically irreducible cubic curve or a cubic surface over a field $K$, with the standard geometric collinearity
relation (1.3) for curves, (3.2.1a) for surfaces, and the binary composition law
(1.1) for curves. For surfaces, we will state the following fancy definition.

\medskip

 {\bf A.1. Definition.} Let $S\subset V(K)$,  and $X_1, \dots ,X_N,\dots $ free commuting but nonassociative variables,
 $$
 w=      (\dots (X_{i_1}\circ X_{i_2})\circ (X_{i_3}\circ  \dots(  \dots\circ  X_{i_k})\dots )
 $$
 a finite word in this variables, 
 $$
 ev\,:\ \{X_1,\dots ,X_N,\dots \}\to S
 $$
 an evaluation map. 
 
 \bigskip
 
 a) A point $p\in V(K)$ is called {\it the strong value of $w$ at $(S, ev )$} if during the inductive calculation
 of 
 $$
 p= ev (W):= (\dots (ev(X_{i_1})\circ ev(X_{i_2}))\circ (ev(X_{i_3}))\circ  \dots (\dots  \circ ev(X_{i_k})\dots )
 $$ 
  we never land in a situation where the result of composition is not uniquely defined, that is $x\circ x$ with singular $x$
  for a curve, or $x\circ y$ where $y=x$ or the line through $x,y$ lies in $V$ for a surface.
  
\bigskip
  
b) A point $p\in S(K)$ is called {\it a weak value of $w$ at $(S, ev )$} if during the inductive calculation
 of 
 $$
 p:= ev_{weak} (W)=(\dots (ev(X_{i_1})\circ ev(X_{i_2}))\circ (ev(X_{i_3}))\circ  \dots (\dots  \circ ev(X_{i_k})\dots  )
 $$ 
 whenever  we  land in a situation  where $\circ$ is not defined,  we are allowed
 to choose as a value of  $y\circ z$ (resp. $y\circ y$) any point of the line $\overline{yz}$ (resp. 
 any point of intersection of a tangent line to $V$ at $y$ with $V$.)
 
Thus, weak evaluation produces a whole set of answers.

\medskip

{\bf A.2. Definition.} {\it (i) A subset  $S\subset V(K)$ {\it strongly generates} $V(K)$,
if $V(K)$ coincides with the set of all strong $S$--values of  all words $w$ as above.

\medskip

(i) A subset  $S\subset V(K)$ {\it weakly generates} $V(K)$,
if $V(K)$ coincides with the set of all weak $S$--values of  all words $w$.}

\medskip

Now we can state two versions of Mordell--Weil problem for cubic surfaces.

\medskip
 {\bf Strong Mordell--Weil problem for $V$:} {\it Is there a finite $S$ that
strongly generates $V(K)$?}

\medskip

{\bf Weak Mordell--Weil problem for $V$:} {\it Is there a finite $S$ that
weakly generates $V(K)$?}

\medskip

For curves, one often calls the weak Mordell--Weil theorem the statement 
that $C(K)/2C(K)$ is finite (referring to the group structure $p+q=e\circ (p\circ q)$).

\medskip

{\bf A.3. Proving strong Mordell--Weil for smooth cubic curves over number fields.}
The classical strategy of proof includes two ingredients.

\smallskip

(a) Introduce an arithmetic {\it height function} $h:\, \P^2(K)\to \bold{R}$.
E.g. for 
$$
K=\Q,\ p=(x_0:x_1:x_2)\in \P^2(\Q),\
x_i\in \bold{Z},\  \roman{g.c.d.}(x_i)=1
$$
put
$$
h(p):= \roman{max}_i |x_i| .
$$

\medskip

(b) Prove {\it the descent property}: 
$\exists H_0\ \roman{such\ that\
if}\ h(p)>H_0,\ 
p\in C(K),\  \roman{then}$
$p=q\circ r\  \roman{for\ some}\ q,r\in C(K)\
\roman{with}\ h(q), h(r)< h(p).
$
\smallskip
The same strategy works for general finitely generated fields. For larger fields,
the strong Mordell--Weil generally fails, but the weak one might  survive.

\medskip

Let $K$ be algebraically closed, or $\bold{R}$, or a finite extension of $\bold{Q}_p$.
Let  $C$ be a smooth cubic curve, $V$ a smooth cubic surface over $K$.

\smallskip

Then:

\medskip

-- $V(K)$ is weakly finitely generated, but not strongly.

\medskip

-- $C(K)$, if non--empty, is not finitely generated.

\medskip

{\bf A.4. Point count on cubic curves.}  It is well known that
as $H\to \infty$,
$$
\roman{card}\,\{p\in C(K)\,|\,h(p)\le H\} = const\cdot (\roman{log}\,H)^{r/2}(1+o(1)),
$$
$$
r:=\roman{rk}\,C(K) =\roman{rk\,Pic}\,C.
\eqno(A.1)
$$

\bigskip

{\bf  A.5. Point count on cubic surfaces.} Here we have only a conjecture and some partial approximations to it:

\medskip
{\bf Conjecture:} as $H\to \infty$,
$$
\roman{card}\,\{p\in V_0(K)\,|\,h(p)\le H\} = const\cdot H(\roman{log}\,H)^{r-1}(1+o(1)),
$$
$$
V_0:= V\setminus \{all \ lines\},\  r:=\roman{rk}\,\roman{Pic}\,V 
\eqno(A.2)
$$
\bigskip

A proof of (A.1) can be obtained by a slight strengthening of the technique
used in the finite generation proof. Namely, one shows that $\roman{log}\,h(p)$
is ``almost a quadratic form'' on $C(K)$. In fact, it differs from a positive defined quadratic
form  by $O(1)$, so that (A.1) follows from the count of lattice point
in an ellipsoid.

\smallskip

The descent property used for Mordell--Weil ensures that this quadratic form is positive definite.

\medskip

How could one attack this conjecture?
For the circle method, there are too few variables.
Moreover, connections with Mordell--Weil for cubic surfaces are totally
missing.

\medskip

Nevertheless, the  inequality
$$
\roman{card}\,\{p\in V_0(K)\,|\,h(p)\le H\} > const\cdot H(\roman{log}\,H)^{r-1}
$$
is proved in [SlSw--D] for cubic surfaces
over $\bold{Q}$ with two rational skew lines. There are also results
for singular surfaces: cf. [Br], [BrD1], [BrD2].

\medskip

{\bf A.6. Some numerical data.} Here I will survey some numerical evidence
computed by Bogdan Vioreanu, cf. [Vi].

\smallskip

In the following tables,  the following notation is used.

\smallskip
{\bf Input/table head:} code $\bold{[a_1,a_2,a_3,a_4]}$ of the surface 
$$V:\,\sum_{i=1}^4 a_ix_i^3=0.$$

\medskip

{\bf Outputs:}

\medskip

(i) $\bold{GEN}$: Conjectural list of weak generators 
$$p:=(x_1:x_2:x_3:x_4)\in V(\Q )$$.

\smallskip

(ii) $\bold{Nr}$: The length of the list $List_{good}$ of all points $x$, ordered by the increasing height
$h(p):=\sum_i |x_i|$, such that any point of the height $\le$ maximal height in $List_{good}$,
is weakly generated by $\bold{GEN}$.

\smallskip

(iii) $\bold{H}_{bad}$: the height of the first point that was NOT shown to be generated by 
$\bold{GEN}$.

\smallskip

(iv) $\bold{L}$: the maximal length of a non--associative word with generators in $(\bold{GEN},\circ )$
one of whose weak values produced an entry in $List_{good}$. 

\medskip

{\bf Example:} For $V=\bold{[1,2,3,4]}$, we have:

\smallskip

$\bold{GEN}=\{p^0:=(1:-1:-1:1)\}$

\smallskip

$\bold{Nr} = 8521$: the first 8521 points in the list of points of increasing height are
weakly generated by the single point $p^0$.

\smallskip

$\bold{L}=13$: the maximal length of a non--associative word in $(p^0, \circ )$
representing some point of the $List_{good}$ was 13.

\smallskip

$\bold{H}_{bad} = 24677$: the first point that was not found to be generated
by $p^0$ was of height 24677.

\bigskip

\centerline{\it SELECTED DATA}

\bigskip

\centerline{$\bold{[1,2,3,5]}, \ \roman{rk\,Pic}=1$}

\centerline{\bf ---------------------------}

\bigskip 

$\bold{GEN}\ \ \ \ \ \ \ \ \ \ \ \ \ \ \ \ \ \ \ \ \ \ \bold{Nr}\ \ \ \ \ \ \ \ \ \ \ \ \ \bold{L}\  \ \ \ \ \ \ \ \ \ \ \ \ \ \bold{H}_{bad}$

{\bf -----------------------------------------------------------------}

(0:1:1:-1) \ \ \ \ \ \ \ \ \ \ \ \ \ \ \ \  15222\ \ \ \ \ \ \   \ \ \  12\ \ \ \ \ \ \ \ \ \ \  \ \ 23243

(1:1:-1:0)

(2:-2:1:1)

\bigskip

\centerline{$\bold{[1,1,5,25]}, \ \roman{rk\,Pic}=2$}

\centerline{\bf ---------------------------}

\bigskip 
 
(1:-1:0:0) \ \ \ \ \ \ \ \ \ \ \ \ \ \ \ \  32419\ \ \ \ \ \ \   \ \ \  9\ \ \ \ \ \ \ \ \ \ \  \ \ \ 30072

(1:4:-2:-1)

\bigskip

\centerline{$\bold{[1,1,7,7]}, \ \roman{rk\,Pic}=2$}

\centerline{\bf ---------------------------}

\bigskip 
 
(0:0:1:-1) \ \ \ \ \ \ \ \ \ \ \ \ \ \ \ \  16063\ \ \ \ \ \ \   \ \ \  7\ \ \ \ \ \ \ \ \ \ \  \ \ \ 2578

(1:-1:0:0)

(1:-1:-1:1)

(1:-1:1:-1)

\bigskip

{\bf A.7. Discussion of other numerical data.} Bogdan Vioreanu studied all in all 16 diagonal  cubic surfaces $V$;
he compiled  lists of all points up to height $10^5$, for some of them
up to height $3\cdot 10^5$.
\smallskip

The conjectural asymptotics (A.2)  seems to be confirmed.

\smallskip

There is a {\it good conjectural expression} for  the constant in (A.2) (for appropriately normalized
height, not the naive one we used). It goes back to works of E.~Peyre.
Its theoretical structure very much reminds the Birch--Swinnerton--Dyer constant for elliptic curves.
For theory and numerical evidence, see [PeT1], [PeT2], [Sw--D2], [Ch-L].

\smallskip
The (weak) finite generation looks confirmed for most of the
considered surfaces, but some stubbornly resist, most notably
$\bold{[17,18,19,20], [4,5,6,7],}$ $\bold{[9,10,11,12]}$.

\smallskip

If one is willing to believe in  weak finite generation (as I am),
the reason for failure might be the following observable fact:

\medskip

When one manages to represent a ``bad" point $p$ of large height
as a non--associative word in the generators $(\bold{GEN},\circ )$,
the height of intermediate results (represented by subwords) tends to be
much higher than $h(p)$, and hence outside of the compiled list of points.

\bigskip

Finally, the relative density of points $p$ for which  ``one--step descent'' works,
that is,
$$
p=q\circ r,\quad h(q), h(r)< h(p),
$$
seems to tend to a certain value $0<d(V)<1.$

\medskip

{\bf Question:} {\it Can one guess a theoretical expression for $d(V)$?}

\medskip

Notice that on each smooth cubic curve $C$, the ``one--step descent'' works for all
points of sufficiently large height, so that $d(C)=1.$

\bigskip

\centerline{\bf References}

\medskip

[Br] T.~D.~Browning. {\it The Manin conjecture in dimension 2.} 

{arXiv:0704.1217}

\smallskip
[BrD1] T.~D.~Browning, U.~ Derenthal. {\it Manin's conjecture for a quartic del Pezzo surface with $A_4$ singularity.} Ann. Inst. Fourier (Grenoble)  59  (2009),  no. 3, 1231--1265. 

\smallskip

[BrD2] T.~D.~Browning, U.~ Derenthal. {\it Manin's conjecture for a cubic surface with $D_5$ singularity.} Int. Math. Res. Not. IMRN 2009, no. 14, 2620--2647. 

\smallskip

[Ch-L] A.~Chambert-Loir. {\it Lectures on height zeta functions:
At the confluence of algebraic geometry, algebraic
number theory, and analysis.} 

arXiv:0812.0947

\smallskip

[HrZ] E.~Hrushovski, B.~Zilber. {\it Zariski geometries.}
Journ. AMS, 9:1 (1996), 1--56.

\smallskip

[K1] D.~S.~Kanevski. {\it Structure of groups, related to cubic surfaces},
Mat. Sb. 103:2, (1977), 292--308 (in Russian); English. transl. in Mat. USSR
Sbornik, Vol. 32:2 (1977), 252--264.

\smallskip

[K2] D.~S.~Kanevsky, {\it On cubic planes and groups connected with cubic
surfaces}. J. Algebra 80:2 (1983), 559--565.

\smallskip

[KaM] D.~S.~Kanevsky, Yu.~I.~Manin.  {\it Composition of points and Mordell--Weil problem for
cubic surfaces. }
In: Rational Points on Algebraic Varieties
(ed. by E.~Peyre, Yu.~Tschinkel), Progress in Mathematics, vol. 199,
Birkh\"auser, Basel, 2001, 199--219. Preprint math.AG/0011198

\smallskip

[M1] Yu.~I.~Manin. {\it Cubic Forms: Algebra, Geometry, Arithmetic.} North
Holland, 1974 and 1986.

\smallskip

[M2] Yu.~I.~Manin. {\it On some groups related to cubic surfaces}. In:
Algebraic Geometry. Tata Press, Bombay, 1968, 255--263.

\smallskip

[M3] Yu.~I.~Manin. {\it Mordell--Weil problem for cubic surfaces}. 
In: Advances in the 
Mathematical Sciences---CRM's 25 Years (L. Vinet, ed.)  CRM Proc. and Lecture Notes, vol. 11, Amer.~Math.~Soc., 
Providence, RI, 1997, pp. 313--318.

\smallskip

[PeT1] E.~Peyre, Yu.~Tschinkel. {\it Tamagawa numbers of diagonal cubic
surfaces, numerical evidence.}  In: Rational Points on Algebraic Varieties,
Progr. Math., 199. Birkh\"auser, Basel, 2001, 275--305. arXiv:9809054

\smallskip

[PeT2] E.~Peyre, Yu.~Tschinkel. {\it Tamagawa numbers of diagonal cubic
surfaces of higher rank.} arXiv:0009092

\smallskip  
[Pr] S.~J.~Pride. {\it Involutary presentations, with applications to
Coxeter groups, NEC-Groups, and groups of Kanevsky}. J. of Algebra 120 (1989),
200--223.

\smallskip

[SlSw--D] J.Slater, H.~P.~F.~Swinnerton--Dyer. {\it Counting points on cubic
surfaces I.} In: Ast\'erisque 251 (1998), 1--12.

\smallskip

[Sw--D1] H.~P.~F.~Swinnerton--Dyer. {\it Universal equivalence for cubic surfaces over
finite and local fields.} Symp.~Math., Bologna 24 (1981), 111--143.

\smallskip

[Sw--D2] H.~P.~F.~Swinnerton--Dyer. {\it Counting points on cubic surfaces II.} In: Geometric 
Methods in Algebra and Number Theory. Progr. Math.,
235, Birkh\"auser, Boston, 2005, pp. 303--309.

\smallskip

[Sw--D3] H.~P.~F.~Swinnerton--Dyer. {\it Universal equivalence for cubic surfaces over
finite and local fields.} Symp.~Math., Bologna 24 (1981), 111--143.

\smallskip

[Vi] B.~G.~Vioreanu. {\it Mordell--Weil problem for cubic surfaces, numerical evidence.}
arXiv:0802.0742

\smallskip

[Z] B.~I.~Zilber. {\it Algebraic geometry via model theory.} Contemp.~Math., vol. 131, Part 3 (1992), 
523--537.

\enddocument